\documentclass[12pt]{amsart}
\usepackage{amsmath,amssymb,latexsym, amsthm, amscd, mathrsfs, stmaryrd} 

\usepackage[all]{xy}
\usepackage{hyperref}
\usepackage{color}
\hypersetup{colorlinks,linkcolor=blue,urlcolor=cyan,citecolor=blue}

\allowdisplaybreaks
\newcommand{\nc}{\newcommand}
\nc{\browntext}[1]{\textcolor{brown}{#1}}
\nc{\greentext}[1]{\textcolor{green}{#1}}
\nc{\redtext}[1]{\textcolor{red}{#1}}
\nc{\bluetext}[1]{\textcolor{blue}{#1}}
\nc{\brown}[1]{\browntext{ #1}}
\nc{\green}[1]{\greentext{ #1}}
\nc{\red}[1]{\redtext{ #1}}
\nc{\blue}[1]{\bluetext{ #1}}
\nc{\zb}[1]{\redtext{From zb: #1}}

\newcommand{\ff}{B}
\newcommand{\tK}{\widetilde{K}}

\newtheorem{thm}{Theorem}  [section]
\newtheorem{cor}[thm]{Corollary}
\newtheorem{lem}[thm]{Lemma}
\newtheorem{prop}[thm]{Proposition}
\newtheorem{conj}[thm]{Conjecture} 

\theoremstyle{remark}
\newtheorem{rem}[thm]{Remark}

\numberwithin{equation}{section}

\newcommand{\mbf}{\mathbf}

\newcommand{\mrm}{\mathrm}

\newcommand{\ev}{\bar{0}}
\newcommand{\odd}{\bar{1}}

\newcommand{\tU}{\widetilde{\mathbf{U}}}

\newcommand{\diag}{\mrm{diag}}

\newcommand{\Aut}{\mrm{Aut}}
\newcommand{\Iblack}{\I_{\bullet}}
\newcommand{\wb}{w_\bullet}

\newcommand{\btau}{\tau}
\def \ty{{\widetilde{y}}}
\def \tk{{\widetilde{k}}}

\newcommand{\N}{\mathbb N}
\newcommand{\K}{\mathbb K}

\newcommand{\Iwhite}{\I_{\circ}}

\newcommand{\ov}{\overline}
\newcommand{\qbinom}[2]{\begin{bmatrix} #1\\#2 \end{bmatrix} }

\newcommand{\U}{\mbf U}

\newcommand{\bs}{{\mbf s}}

\newcommand{\Ui}{{\mbf U}^\imath}

\newcommand{\vs}{\varsigma}

\newcommand{\Z}{\mathbb Z}
\newcommand{\I}{ I}
\def\Id{\mathrm{Id}}

\newcommand{\T}{\mbf T}

\def \fg{\mathfrak{g}}

\nc{\fprime}{\bold{'f}}

\def \bU{{\mathbf U}}
\def \bvs{{\boldsymbol{\varsigma}}}

\newcommand{\tUi}{\widetilde{{\mathbf U}}^\imath}

\begin{document}

\title[Serre-Lusztig relations for  $\imath$quantum groups III]{Serre-Lusztig relations for  $\imath$quantum groups III}

\author[Xinhong Chen]{Xinhong Chen}
\address{Department of Mathematics, Southwest Jiaotong University, Chengdu 610031, P.R.China}
\email{chenxinhong@swjtu.edu.cn}

\author[Ming Lu]{Ming Lu}
\address{Department of Mathematics, Sichuan University, Chengdu 610064, P.R.China}
\email{luming@scu.edu.cn}

\author[Weiqiang Wang]{Weiqiang Wang}
\address{Department of Mathematics, University of Virginia, Charlottesville, VA 22904}
\email{ww9c@virginia.edu}



\begin{abstract}
let $\widetilde{\bf U}^\imath$ be a quasi-split universal $\imath$quantum group associated to a quantum symmetric pair $(\widetilde{\bf U}, \widetilde{\bf U}^\imath)$ of Kac-Moody type with a diagram involution $\tau$. We establish the Serre-Lusztig relations for $\widetilde{\bf U}^\imath$ associated to a simple root $i$ such that $i \neq \tau i$, complementary to the Serre-Lusztig relations associated to $i=\tau i$ which we obtained earlier. A conjecture on braid group symmetries on $\widetilde{\bf U}^\imath$ associated to $i$ disjoint from $\tau i$ is formulated.
\end{abstract}

\maketitle

\section{Introduction}

\subsection{}

Lusztig \cite[Chapter 7]{Lus93} formulated the higher order quantum Serre relations for Drinfeld-Jimbo quantum groups, which we shall refer to as Serre-Lusztig relations. They are intimately related to the braid group actions on quantum groups. 

Associated to a Satake diagram $(I =I_\circ \cup I_\bullet, \tau)$ where $\tau$ is a diagram involution, a quantum symmetric pair $(\U, \Ui)$ \cite{Le99, Le03} consists of a Drinfeld-Jimbo quantum group $\U$ and its coideal subalgebra $\Ui$. We refer to $\Ui$ as an $\imath$quantum group, and call it {\em quasi-split} if $I_\bullet =\emptyset$. A universal $\imath$quantum group $\tUi$ introduced in \cite{LW19a} has Chevalley generators $B_i, \tk_i$ ($i \in I$), and Letzter's $\imath$quantum groups with parameters are obtained from universal $\imath$quantum groups by central reductions. 

Generalizing \cite{Le03, Ko14, BK15}, a Serre presentation for quasi-split $\imath$quantum groups of arbitrary Kac-Moody type was obtained in \cite{CLW18}, where the $\imath$Serre relations are expressed in terms of the $\imath$divided powers \cite{BW18a, BeW18}.
The authors formulated in \cite{CLW21} the Serre-Lusztig relations between $B_i, B_j$ for $i\neq j \in I$, for (mostly) quasi-split universal $\imath$quantum groups in the Kac-Moody generality, for $i =\tau i$; also see \cite{BV15} for some earlier attempt and examples. These are higher order relations associated to the $\imath$Serre relations given in \cite{CLW18}, where the $\imath$divided powers (associated to $i =\tau i$) are again indispensible. These Serre-Lusztig relations are more involved in both formulations and proofs than their quantum group counterparts. A further generalization of the Serre-Lusztig relations for $i=\tau i$ was obtained in \cite{CLLW22}, where the ``quasi-split" condition on $\imath$quantum groups was completely removed. 


%
%
\subsection{}

In this paper, we shall establish the Serre-Lusztig relations between $B_i, B_j$, for $i\neq j \in I$, for (mostly) quasi-split universal $\imath$quantum groups of Kac-Moody type, for $i \neq \tau i$; this complements the Serre-Lusztig relations for $\widetilde{\bf U}^\imath$ with $i=\tau i$ in \cite{CLW21, CLLW22}. We then formulate a conjecture on closed formulas for the braid group symmetries associated with $i\neq \tau i$ on $\tUi$. All the formulas are expressed in terms of the $\imath$divided powers. 
\subsection{}

It turns out that the Serre-Lusztig relations between $B_i, B_j$, for $\tau i \neq i\neq j \in I$, are standard as for quantum groups  {\em unless $j =\tau i$}; see Proposition~\ref{prop:LS2}. 

So we restrict our discussion to the relations between $B_i, B_{\tau i}$ ($i\neq \tau i$). 
A Serre type relation in $\tUi$ between $B_i, B_{\tau i}$ which contains a term involving $\tk_i$ and a second term involving $\tk_{\tau i}$ is given in \eqref{relation5}. This relation is a universal variant of a relation established in \cite[Theorem~3.6]{BK15} for $\imath$quantum groups $\Ui$ of arbitrary Kac-Moody type, generalizing a relation due to Letzter \cite[Theorem~7.1]{Le03} in finite type. We shall refer to \eqref{relation5} as the BKL relation. The BKL relation has the minimal degree in the hierarchy of Serre-Lusztig relations between $B_i, B_{\tau i}$, which we are going to formulate.

Denote by $(c_{ij})_{i,j \in I}$ the generalized Cartan matrix. 
We introduce a family of elements in $\tUi$, denoted by $\ty_{i,\tau i;1,m,e}$ in \eqref{eq:fitaui}, for $m\ge 0$ and $e =\pm 1$. In this notation, the BKL relation becomes $\ty_{i,\tau i;1,1-c_{i,\tau i},1} =0$. We establish a recursive relation for $\ty_{i,\tau i;1,m,e}$ ($m\ge 0$) in Theorem~\ref{thm:recursion}, from which we derive the Serre-Lusztig relations 
\[
\ty_{i,\tau i;1,m,e} =0, \qquad \text{ for } m\ge 1-c_{i,\tau i}.
\] 
The elements $\ty_{i,\tau i;1,m,e}$ admit a simpler reformulation, for $m > 1-c_{i,\tau i}$, in the sense that it does not contain a term involving $\tk_{ i}$ for $e=1$ (respectively, $\tk_{\tau i}$ for $e=-1$). In this way, the Serre-Lusztig relations for $m > 1-c_{i,\tau i}$ look strikingly different from the BKL relation (= Serre-Lusztig for $m =1-c_{i,\tau i}$); see Theorem~\ref{thm:HOSII}.

\subsection{}

Braid group symmetries \cite[Part~V]{Lus93} are fundamental in Drinfeld-Jimbo quantum groups. Formulas for braid group actions on (mostly) quasi-split $\imath$quantum groups of finite type (in distinguished parameters) have been obtained in \cite{KP11} via computer computations. A conceptual approach via reflection functors for braid group symmetries on quasi-split $\tUi$ of finite type was developed in \cite{LW21a}. Dobson \cite{D20} obtained braid group actions for $\imath$quantum groups of type AIII, which are not necessarily quasi-split. It is a basic problem to formulate conceptually the braid group symmetries on $\imath$quantum groups in great generalities; this has been achieved in \cite{WZ22} for (not necessarily quasi-split) $\imath$quantum groups of arbitrary finite type. 

The Serre-Lusztig relations for $i=\tau i$ have led to a conjecture \cite[Conjecture~ 6.5]{CLW21}  on braid group operators $\T''_{i,e}$ and $\T'_{i,e}$, for $i=\tau i$ and $e=\pm 1$. 
In contrast, for weight reason (see Remark~\ref{rem:weight}), the Serre-Lusztig relations obtained above do not allow us to guess directly formulas for $\T''_{i,e}$ and $\T'_{i,e}$, for $i\neq \tau i$ and $e=\pm 1$. From a general consideration of restricted Weyl groups, we only expect the braid operators for $i\neq \tau i$ exist when the Cartan integers $c_{i, \tau i} =0, -1$; cf. \cite{KP11, Lus03}, and we know what the weight $\T''_{i,e}(B_j)$ should be. With considerable efforts, in case $c_{i, \tau i} =0$ we formulate a conjecture on closed formulas of the automorphisms $\T''_{i,e}$ and $\T'_{i,e}$, for $i=\tau i$ and $e=\pm 1$; see Conjecture~\ref{conj:BG1}. We observe a feature similar to the Serre-Lusztig relations that in any given braid group formula acting on $B_j$ only terms involving $\tk_{ i}$ or $\tk_{\tau i}$ (but not both) show up.

We remark that the Serre-Lusztig relations and the conjecture on the braid group automorphisms in this paper are valid for some more general (beyond quasi-split) $\imath$quantum groups; see \S\ref{subsec:nonQS}. We choose to formulate the main body of this paper in the quasi-split setting in which the notations are simplified and the results are most complete. \cite[Conjecture~ 6.5]{CLW21} and (a major portion of) Conjecture~\ref{conj:BG1} in this paper will be established in a forthcoming work \cite{LW21b} for {\em quasi-split} $\imath$quantum groups (under some additional assumptions on Cartan integers; see Remark~\ref{rem:weight2}) via a Hall algebra approach. Having the exact (conjectural) formulas available is crucial for the Hall algebra approach to work. To establish the conjectures in full generality shall motivate some future developments.

\subsection{}
This paper is organized as follows.
In Section~\ref{sec:prelim}, we review and set up notations for quantum groups and $\imath$quantum groups. In Section~\ref{sec:SL}, we establish the Serre-Lusztig relations for $\imath$quantum groups and formulate a conjecture on braid group symmetries, when $i \neq \tau i$.

\vspace{2mm}
{\bf Acknowledgement.}
 ML thanks Shanghai Key Laboratory of Pure Mathematics and
Mathematical Practice, East China Normal University for hospitality and support.
XC is partially supported by the Fundamental Research Funds for the Central Universities grant No. 2682020ZT100. ML is partially supported by the Science and Technology Commission of Shanghai Municipality (grant No. 18dz2271000), and the National Natural Science Foundation of China (grant No. 12171333). WW is partially supported by the NSF grant DMS-2001351.

\section{Quantum symmetric pairs and $\imath$quantum groups}
  \label{sec:prelim}

In this section, we recall the definitions of $\imath$quantum groups and universal $\imath$quantum groups arising from quantum symmetric pairs. 

\subsection{Quantum groups}

Let $C=(c_{ij})_{i,j \in I}$ be a symmetrizable generalized Cartan matrix with its symmetrizer $D=\diag(\epsilon_i\mid \epsilon_i\in \Z_{\ge 1},\; i\in I)$, i.e., $DC$ is symmetric.
Let $\fg$ be the corresponding Kac-Moody Lie algebra. Let $\alpha_i$ ($i\in I $) be the simple roots of $\fg$, and denote the root lattice by $\Z^\I:=\Z\alpha_1\oplus\cdots\oplus\Z\alpha_n$. The {\em simple reflection} $s_i:\Z^{ I}\rightarrow\Z^{ I}$ is defined to be $s_i(\alpha_j)=\alpha_j-c_{ij}\alpha_i$, for $i,j\in  I$.
Denote the Weyl group by $W =\langle s_i\mid i\in  I\rangle$.

Let $q$ be an indeterminate, and denote
\[
q_i:=q^{\epsilon_i}, \qquad \forall i\in I.
\]
For $n,m\in \Z$ with $m\ge 0$ and indeterminate $t$, we denote the quantum integers and quantum binomial coefficients as
\begin{align*}
\begin{split}
[n]_t =\frac{t^n-t^{-n}}{t-t^{-1}},
\quad
[m]_t^! =\prod_{i=1}^m [i]_t,
\quad
\qbinom{n}{d}_t  &=
\begin{cases}
\frac{[n]_t[n-1]_t \ldots [n-d+1]_t}{[d]_t^!}, & \text{ if }d \ge 0,
\\
0, & \text{ if }d<0.
\end{cases}
\end{split}
\end{align*}
We shall use these notations by setting $t=q$ or $q_i$. 

Let $\K$ be a field of characteristic $0$. Assume that a symmetrizable generalized Cartan matrix $C$ is given. 
Then $\tU := \tU_q(\fg)$ is the associative $\K(q)$-algebra with generators $E_i,F_i,\tK_i,\tK_i'$ for all $i\in I$ where $\tK_i,\tK'_i$ are invertible, subject to the following relations:
\begin{align}
&\tK_i\tK_{j}=\tK_j\tK_{i},\; \tK_i\tK'_{j}=\tK'_{j}\tK_i,\; \tK'_i\tK'_{j}=\tK'_j\tK'_{i}, \quad \forall i,j\in I,\label{Q1}
\\
&\tK_i E_j=q_i^{c_{ij}} E_j \tK_i,  \qquad \tK_i F_j=q_i^{-c_{ij}} F_j \tK_i,
\label{eq:EK}
\\
&\tK_i' E_j=q_i^{-c_{ij}} E_j \tK_i',  \qquad \tK_i' F_j=q_i^{c_{ij}} F_j \tK_i',
 \label{eq:K2}
 \\
 &[E_i,F_j]=\delta_{ij}\frac{\tK_{i}-\tK_{i}'}{q_i-q_i^{-1}},
\label{Q4}\\
&  
\sum_{n=0}^{1-c_{ij}} (-1)^n  E_i^{(n)} E_j E_i^{({1-c_{ij}-n})}=
\sum_{n=0}^{1-c_{ij}} (-1)^n  F_i^{(n)}F_j F_i^{(1-c_{ij}-n)}=0,
 \quad \forall i \neq j\in I.
\label{q serre}
\end{align}
In the above $q$-Serre relations, we have used the divided powers
\[
F_i^{(n)} =F_i^n/[n]_{q_i}^!, \quad E_i^{(n)} =E_i^n/[n]_{q_i}^!, \quad \text{ for } n\ge 1 \text{ and  } i\in I.
\]

Note that $\tK_i\tK'_i$ are central in $\tU$ for any $i\in I$.  The comultiplication $\Delta: \widetilde{\U} \rightarrow \widetilde{\U} \otimes \widetilde{\U}$ is defined as follows:
\begin{align}  \label{eq:Delta}
\begin{split}
\Delta(E_i)  = E_i \otimes 1 + \tK_i \otimes E_i, & \quad \Delta(F_i) = 1 \otimes F_i + F_i \otimes \tK_{i}', \\
 \Delta(\tK_{i}) = \tK_{i} \otimes \tK_{i}, & \quad \Delta(\tK_{i}') = \tK_{i}' \otimes \tK_{i}'.
 \end{split}
\end{align}
The Chevalley involution $\omega$ on $\tU$ is given by
\begin{align}  \label{eq:omega}
\omega(E_i)  = F_i,\quad  \omega(F_i) = E_i,\quad \omega(\tK_{i}) = \tK_{i}' , \quad \omega(\tK_{i}') =\tK_{i}, \quad \forall i\in I.
\end{align}

Analogously as for $\tU$, the quantum group $\bU$ is defined to be the $\K(q)$-algebra generated by $E_i,F_i, K_i, K_i^{-1}$, for all $i\in I$, subject to the  relations modified from \eqref{Q1}--\eqref{Q4} with $\tK_i$ and $\tK_i'$ replaced by $K_i$ and $K_i^{-1}$, respectively. The comultiplication $\Delta$ and Chevalley involution $\omega$ on $\U$ are obtained by modifying \eqref{eq:Delta}--\eqref{eq:omega} with $\tK_i$ and $\tK_i'$ replaced by $K_i$ and $K_i^{-1}$, respectively (cf. \cite{Lus93}; beware that our $\tK_i$ has a different meaning from $\tK_i \in \U$ therein.) 
We have ${\bU} \cong \widetilde{\bU} \big/ (\tK_i \tK_i' -1 \mid   i\in I)$. 


%
\subsection{The $\imath$quantum groups $\tUi$ and $\Ui$}
  \label{subsec:iQG}

For a  (generalized) Cartan matrix $C=(c_{ij})$, let $\btau$ be an involution in $\Aut(C)$, i.e., a permutation of $I$ such that $c_{ij}=c_{\btau i,\btau j}$, and $\btau^2=\Id$. 
We define $\widetilde{\bU}^\imath$ to be the $\K(q)$-subalgebra of $\tU$ generated by
\[
B_i= F_i +  E_{\btau i} \tK_i',
\qquad \tk_i= \tK_i \tK_{\btau i}', \quad \forall i \in I.
\]
According to \cite{LW19a}, the elements $\tk_i$ (for $\btau i=i$) and $\tk_i \tk_{\btau i}$  (for $i\neq \btau i$) are central in $\tUi$.

Let $\bvs=(\vs_i)\in  (\K(q)^\times)^{I}$ be such that $\vs_i=\vs_{\btau i}$, for each $i\in I$ which satisfies $c_{i, \btau i}=0$.
Let $\Ui:=\Ui_{\bvs}$ be the $\K(q)$-subalgebra of $\bU$ generated by
\[
B_i= F_i+\vs_i E_{\btau i}K_i^{-1},
\quad
k_j= K_jK_{\btau j}^{-1},
\qquad  \forall i \in I,j\in I  \text{ such that }\btau j\neq j.
\]
It is known \cite{Le99, Ko14} that $\bU^\imath$ is a right coideal subalgebra of $\bU$ in the sense that $\Delta: \Ui \rightarrow \Ui\otimes \U$; and $(\bU,\Ui)$ is called a \emph{quantum symmetric pair}, as they specialize at $q=1$ to $(\U(\fg), \U(\fg^{\omega \btau}))$, where $\btau$ is understood here as an automorphism of $\fg$. Let
\begin{align}
\label{eq:ci}
I_\tau = \{ \text{fixed representatives of $\tau$-orbits in $I$} \}.
\end{align}

The algebra $\widetilde{\bU}^\imath$ is a right coideal subalgebra of $\widetilde{\bU}$ \cite{LW19a}. 
The algebras $\Ui_{\bvs}$, for $\bvs \in  (\K(q)^\times)^{I}$, are obtained from $\tUi$ by central reductions as follows. 

\begin{prop}
	[\text{\cite[Proposition 6.2]{LW19a}}]
 \label{prop:QSP12}
The algebra $\Ui$ is isomorphic to the quotient of $\tUi$ by the ideal generated by
\begin{align*}
\tk_i - \vs_i \; (\text{for } i =\btau i),
\qquad  \tk_i \tk_{\btau i} - \vs_i \vs_{\btau i}  \;(\text{for } i \neq \btau i).
\end{align*}
The isomorphism is given by sending $B_i \mapsto B_i, k_j \mapsto \frac{1}{\vs_{\btau j}} \tk_j, k_j^{-1} \mapsto \frac{1}{\vs_{ j}} \tk_{\btau j}, \forall i\in I, j\in I\backslash I_\tau$.
\end{prop}

\subsection{A Serre presentation of $\tUi$}\label{subsection:iserre presentation}

For  $i\in I$ with $\btau i\neq i$, imitating Lusztig's divided powers,  we define the {\em divided power} of $B_i$ to be
\begin{align}
  \label{eq:iDP1}
  B_i^{(m)}:=B_i^{m}/[m]_{q_i}^!, \quad \forall m\ge 0, \qquad (\text{if } i \neq \tau i).
\end{align}

For $i\in I$ with $\btau i= i$, generalizing \cite{BeW18}, we define the {\em $\imath$divided powers} of $B_i$ to be
\begin{align}
&&\ff_{i,\odd}^{(m)}=\frac{1}{[m]_{q_i}^!}\left\{ \begin{array}{ccccc} B_i\prod_{j=1}^k (B_i^2-q_i\tk_i[2j-1]_{q_i}^2 ) & \text{if }m=2k+1,\\
\prod_{j=1}^k (B_i^2-q_i\tk_i[2j-1]_{q_i}^2) &\text{if }m=2k; \end{array}\right.
  \label{eq:iDPodd} \\
  \notag\\
&&\ff_{i,\ev}^{(m)}= \frac{1}{[m]_{q_i}^!}\left\{ \begin{array}{ccccc} B_i\prod_{j=1}^k (B_i^2-q_i\tk_i[2j]_{q_i}^2 ) & \text{if }m=2k+1,\\
\prod_{j=1}^{k} (B_i^2-q_i\tk_i[2j-2]_{q_i}^2) &\text{if }m=2k. \end{array}\right.
 \label{eq:iDPev}
\end{align}



Denote
\[
(a;x)_0=1, \qquad (a;x)_n =(1-a)(1-ax)  \cdots (1-ax^{n-1}), \quad \forall n\ge 1.
\]

\begin{prop}[Serre presentation of universal quasi-split $\imath$quantum groups \cite{CLW18}]\label{prop:Serre}
The $\K(q)$-algebra $\tUi$ has a presentation with generators $B_i$, $\tk_i$ $(i\in I)$ and the relations \eqref{relation1}--\eqref{relation6} below: for $\ell \in I$, and $i\neq j \in I$,
\begin{align}
\tk_i \tk_\ell =\tk_\ell \tk_i, \qquad
\tk_\ell B_i&  = q_i^{c_{\btau \ell,i} -c_{\ell i}} B_i \tk_\ell,
\label{relation1}
\\
B_iB_{j}-B_jB_i =0, \quad \text{ if }&c_{ij} =0 \text{ and }\btau i\neq j,\label{relation2}
\\
\sum_{n=0}^{1-c_{ij}} (-1)^nB_i^{(n)}B_jB_i^{(1-c_{ij}-n)} &=0, \qquad \text{ if } j \neq \btau i\neq i, \label{relation3}
\\
\sum_{n=0}^{1-c_{i,\btau i}} (-1)^{n+c_{i,\tau i}}B_i^{(n)}B_{\btau i}B_i^{(1-c_{i,\btau i}-n)}& = 
\label{relation5}
\\
\frac{1}{q_i-q_i^{-1}} \Big(q_i^{c_{i,\btau i}} (q_i^{-2};q_i^{-2})_{-c_{i,\btau i}} 
& B_i^{(-c_{i,\btau i})} \tk_i  
  -(q_i^{2};q_i^{2})_{-c_{i,\btau i}}B_i^{(-c_{i,\tau i})} \tk_{\btau i}  \Big),\quad
\text{ if } \btau i \neq i,
 \notag \\
\sum_{r=0}^{1-c_{ij}} (-1)^r  B_{i,\overline{p_i}}^{(r)}B_j B_{i,\overline{p}_i+\ov{c_{ij}}}^{(1-c_{ij}-r)} &=0,\qquad   \text{ if }\btau i=i.
\label{relation6}
\end{align}
\end{prop}

\begin{rem}
The presentation of $\Ui$ of finite type was due to G. Letzter \cite[Theorem~7.1]{Le03}. 
In the setting of $\imath$quantum group $\Ui$ of Kac-Moody type, the relation \eqref{relation5} was established in \cite[Theorem 3.6]{BK15} (generalizing a relation in finite type \cite{Le03}) and will be referred to as {\em the BKL relation}. This relation in case $c_{i,\tau i}=0$ is essentially the standard quantum $\mathfrak{sl}_2$ relation between $E$ and $F$; see the proof of Corollary~\ref{cor:rank1}. 
\end{rem}

\begin{lem}
\label{lem:bar}
(a) There exists a $\K$-algebra automorphism $ \psi_\imath: \tUi\rightarrow \tUi$ (called a bar involution) such that
\[
\psi_\imath(q)=q^{-1}, \quad
\psi_\imath(\tk_i)=q_i^{c_{i,\tau i}}\tk_{\tau i}, \quad
\psi_\imath (B_i)=B_i,  \quad
\forall i\in I.
\]
(b) 
There exists a $\K(q)$-algebra anti-involution $\sigma:\tUi\rightarrow \tUi$ such that
	$$\sigma(B_i)=B_i, \quad \sigma(\tk_i)= \tk_{\btau i}, 
	\quad \forall i\in I.$$
\end{lem}

\begin{proof}
(a) It suffices to show that $\psi_\imath$ preserves all the defining relations for $\tUi$ in Proposition~\ref{prop:Serre}. 
Note that if $i=\btau i$, then $c_{i,\tau i}=2$, and $\psi_\imath (\tk_i)=q_i^2\tk_i$ in this case. So $\psi_\imath$ fixes $B_{i,\ov{p}}^{(n)}$ in \eqref{eq:iDPodd}--\eqref{eq:iDPev}, for any $i\in I$ so that $i=\btau i$, $\ov{p}\in\Z_2$ and $n\in\N$. Hence, one checks readily that $\psi_\imath$ preserves the relations \eqref{relation1}--\eqref{relation6} except perhaps \eqref{relation5}. For \eqref{relation5}, as its LHS is clearly fixed by $\psi_\imath$, it suffices to check that its RHS is fixed by $\psi_\imath$ as follows:
\begin{align*}
&\psi_\imath \left(\frac{1}{q_i-q_i^{-1}}\left(q_i^{c_{i,\btau i}} (q_i^{-2};q_i^{-2})_{-c_{i,\btau i}} B_i^{(-c_{i,\btau i})} \tk_i \right.
 \left. -(q_i^{2};q_i^{2})_{-c_{i,\btau i}}B_i^{(-c_{i,\tau i})} \tk_{\btau i}  \right)\right)
 \\
&= -\frac{1}{q_i-q_i^{-1}}\Big(q_i^{-c_{i,\btau i}} (q_i^{2};q_i^{2})_{-c_{i,\btau i}} q_i^{c_{i,\btau i}} B_i^{(-c_{i,\btau i})}\tk_{\btau i} - (q_i^{-2};q_i^{-2})_{-c_{i,\btau i}}q_i^{c_{i,\btau i}}B_i^{(-c_{i,\tau i})}\tk_i\Big)
\\
&= \frac{1}{q_i-q_i^{-1}}\Big( q_i^{c_{i,\btau i}}(q_i^{-2};q_i^{-2})_{-c_{i,\btau i}}B_i^{(-c_{i,\tau i})}\tk_i-(q_i^{2};q_i^{2})_{-c_{i,\btau i}} B_i^{(-c_{i,\btau i})}\tk_{\btau i} \Big).
\end{align*}

(b)  It also follows by inspection that $\sigma$ preserves the relations for $\tUi$ in Proposition~\ref{prop:Serre}. The anti-involution $\sigma$ is given by the composition of $\tau$ and the anti-involution $\sigma_\imath$ in \cite[Lemma~ 2.3]{CLW21}. (We thank Weinan Zhang for a helpful discussion.)
\end{proof}

\section{Serre-Lusztig relations for $\tUi$}
\label{sec:SL}

In this section, we shall formulate and establish the Serre-Lusztig relations associated to the BKL relation \eqref{relation5} in $\tUi$, for $i \in I$ such that $\tau i \neq i$. The Serre-Lusztig relations associated to the standard Serre relation \eqref{relation3} is also given. 

\subsection{A recursive formula}

For any $i\in I$ such that $\tau i\neq i$ and $m\in \N$, recalling the divided powers $B_i^{(m)}$ from \eqref{eq:iDP1}, we define
\begin{align}
\label{eq:fitaui}
\ty_{i,\tau i;1,m,e} &:= \sum_{r+s=m} (-1)^{r+c_{i,\tau i}}q_i^{er(1-c_{i,\tau i}-m)}B_i^{(r)}B_{\tau i}B_i^{(s)}- \frac{[1-c_{i,\tau i}]_{q_i}^!}{[m]_{q_i}}(q_i-q_i^{-1})^{-c_{i,\tau i}-1} \times
\\
\notag
&\quad \Big\{ 
\prod_{j=0}^{c_{i,\tau i}+m-2} \big(-q_i^{e(2j-c_{i,\tau i}-2m+2)}+q_i^{c_{i,\tau i}-2}\big)
q_i^{\frac{-c_{i,\tau i}^2+3c_{i,\tau i}}{2}}B_{i}^{(m-1)}\tk_i
\\\notag
&\qquad -(-1)^{c_{i,\tau i}}\prod_{j=0}^{c_{i,\tau i}+m-2} \big(-q_i^{e(2j-c_{i,\tau i}-2m+2)}+q_i^{2-c_{i,\tau i}}\big)q_i^{\frac{c_{i,\tau i}^2-c_{i,\tau i}}{2}}B_{i}^{(m-1)}\tk_{\tau i}\Big\}.
\end{align}
%
It is understood that the above expression $\prod\limits_{j=0}^{c_{i,\tau i}+m-2} (***)=1$ if $c_{i,\tau i}+m-2<0$. In other words, for $m \le 1 -c_{i,\tau i}$, we have 
\begin{align*}
\ty_{i,\tau i;1,m,e} &= \sum_{r+s=m} (-1)^{r+c_{i,\tau i}}q_i^{er(1-c_{i,\tau i}-m)}B_i^{(r)}B_{\tau i}B_i^{(s)}   
\\
\notag
&\quad 
- \frac{[1-c_{i,\tau i}]_{q_i}^!}{[m]_{q_i}}(q_i-q_i^{-1})^{-c_{i,\tau i}-1} \Big\{ 
q_i^{\frac{-c_{i,\tau i}^2+3c_{i,\tau i}}{2}}B_{i}^{(m-1)}\tk_i
 -(-1)^{c_{i,\tau i}} q_i^{\frac{c_{i,\tau i}^2-c_{i,\tau i}}{2}}B_{i}^{(m-1)}\tk_{\tau i}
 \Big\}.
\end{align*}
As we shall see, the BKL relation~\eqref{relation5} can be reformulated as $\ty_{i,\tau i;1,1-c_{i,\tau i},e}=0$. 

For $m > 1 -c_{i,\tau i}$, the $\ty_{i,\tau i;1,m,e}$ can also be much simplified; see Theorem~\ref{thm:HOSII} below.

Recall the anti-involution $\sigma$ of $\tUi$ from Lemma~\ref{lem:bar}. Define
\begin{align}  \label{eq:ys}
\ty'_{i,\tau i;1,m,e}:=\sigma(\ty_{i,\tau i;1,m,e}),
\qquad
\forall \tau i \neq i \in I. 
\end{align}

We have the following recursive relations among $\ty_{i,\tau i;1,m,e}$. Guessing the right definition of $\ty_{i,\tau i;1,m,e}$ is (a most difficult) part of the statement!

\begin{thm}
  \label{thm:recursion}
Let $i\in \I$ be such that  $\btau i\neq i$. Then for any  $m\in \N$, and $e=\pm1$, we have
\begin{align}
	\label{eq:RR1}
-q_i^{-e(2m+c_{i,\tau i})} & B_i \ty_{i,\tau i;1,m,e} +\ty_{i,\tau i;1,m,e}B_i
=[m+1]_{i}\, \ty_{i,\tau i;1,m+1,e},
		\\
	\label{eq:RR2}
-q_i^{-e(2m+c_{i,\tau i})} & \ty'_{i,\tau i;1,m,e}B_i+B_i\ty'_{i,\tau i;1,m,e}
=[m+1]_{i}\, \ty'_{i,\tau i;1,m+1,e}.
\end{align}
\end{thm}

\begin{proof}
The identity \eqref{eq:RR2} follows by applying the anti-involution $\sigma$ to \eqref{eq:RR1}; see \eqref{eq:ys}. Thus, it suffices to prove \eqref{eq:RR1}. 

By a formal computation (cf. \cite[Lemma 7.1.2]{Lus93}), we have
\begin{align*}
-q_i^{e(-c_{i,\tau i}-2m)} & B_i\Big( \sum_{r+s=m} (-1)^{r+c_{i,\tau i}}q_i^{er(1-c_{i,\tau i}-m)}B_i^{(r)}B_{\tau i}B_i^{(s)}\Big)
	\\
&+\Big( \sum_{r+s=m} (-1)^{r+c_{i,\tau i}}q_i^{er(1-c_{i,\tau i}-m)}B_i^{(r)}B_{\tau i}B_i^{(s)}\Big)B_i	\\
&=[m+1]_i \sum_{r+s=m+1} (-1)^{r+c_{i,\tau i}}q_i^{er(-c_{i,\tau i}-m)}B_i^{(r)}B_{\tau i}B_i^{(s)}.
\end{align*}

On the other hand, by \eqref{relation1}, we have
\begin{align*}
\tk_iB_i=q_i^{c_{i,\tau i}-2} B_i \tk_i \text{ and }
\tk_{\tau i}B_i=q_i^{-c_{i,\tau i}+2} B_i \tk_{\tau i}.
\end{align*}
Using these identities and the definition of divided powers \eqref{eq:iDP1}, we verify by a direct computation that  
\begin{align*}
 -q_i^{e(-c_{i,\tau i}-2m)}B_i \cdot 
  & \Big\{\frac{[1-c_{i,\tau i}]_{q_i}^!}{[m]_{q_i}}(q_i-q_i^{-1})^{-c_{i,\tau i}-1}\\
	&\times \Big( \prod_{j=0}^{c_{i,\tau i}+m-2} \big(-q_i^{e(2j-c_{i,\tau i}-2m+2)}+q_i^{c_{i,\tau i}-2}\big)q_i^{\frac{-c_{i,\tau i}^2+3c_{i,\tau i}}{2}}B_{i}^{(m-1)}\tk_i\\
	&-(-1)^{c_{i,\tau i}}\prod_{j=0}^{c_{i,\tau i}+m-2} \big(-q_i^{e(2j-c_{i,\tau i}-2m+2)}+q_i^{2-c_{i,\tau i}}\big)q_i^{\frac{c_{i,\tau i}^2-c_{i,\tau i}}{2}}B_{i}^{(m-1)}\tk_{\tau i}\Big)\Big\}\\
 + \frac{[1-c_{i,\tau i}]_{q_i}^!}{[m]_{q_i}} & (q_i-q_i^{-1})^{-c_{i,\tau i}-1}\\
	&\times \Big( \prod_{j=0}^{c_{i,\tau i}+m-2} \big(-q_i^{e(2j-c_{i,\tau i}-2m+2)}+q_i^{c_{i,\tau i}-2}\big)q_i^{\frac{-c_{i,\tau i}^2+3c_{i,\tau i}}{2}}B_{i}^{(m-1)}\tk_i\\
	&  -(-1)^{c_{i,\tau i}}\prod_{j=0}^{c_{i,\tau i}+m-2} \big(-q_i^{e(2j-c_{i,\tau i}-2m+2)}+q_i^{2-c_{i,\tau i}}\big)q_i^{\frac{c_{i,\tau i}^2-c_{i,\tau i}}{2}}B_{i}^{(m-1)}\tk_{\tau i}\Big) \cdot B_i\\
	=  [m+1]_{q_i} & \frac{[1-c_{i,\tau i}]_{q_i}^!}{[m+1]_{q_i}}(q_i-q_i^{-1})^{-c_{i,\tau i}-1}\\
	&\times \Big\{ \prod_{j=0}^{c_{i,\tau i}+m-1} \big(-q_i^{e(2j-c_{i,\tau i}-2m)}+q_i^{c_{i,\tau i}-2}\big)q_i^{\frac{-c_{i,\tau i}^2+3c_{i,\tau i}}{2}}B_{i}^{(m)}\tk_i\\
	& -(-1)^{c_{i,\tau i}}\prod_{j=0}^{c_{i,\tau i}+m-1} \big(-q_i^{e(2j-c_{i,\tau i}-2m)}+q_i^{2-c_{i,\tau i}}\big)q_i^{\frac{c_{i,\tau i}^2-c_{i,\tau i}}{2}}B_{i}^{(m)}\tk_{\tau i}\Big\}.
\end{align*}

Combining the above two formulas and using \eqref{eq:fitaui}, we have proved the formula \eqref{eq:RR1}.
\end{proof}

\subsection{Serre-Lusztig for $i\neq \tau i$}

It turns out that one of the 2 messy products in the definition \eqref{eq:fitaui} of $\ty$ is 0 when $m$ is large enough (more precisely, when $m>1-c_{i,\tau i}$). 

\begin{thm}
\label{thm:HOSII}
	Let $i\in \I$ be such that  $\btau i\neq i$. 
For $m\geq 1-c_{i,\tau i}$, we have
\begin{align}
	\label{eq:y=0}
\ty_{i,\tau i;1,m,e}=0,\qquad \text{ and } \;\; 
\ty'_{i,\tau i;1,m,e}=0.
\end{align}
Equivalently, for $m>1-c_{i,\tau i}$, we have
 \begin{align}
 \label{eq:HOSII1}
\sum_{r+s=m} & (-1)^{r+c_{i,\tau i}} q_i^{r(1-c_{i,\tau i}-m)}B_i^{(r)}B_{\tau i}B_i^{(s)}\\  
&= (-1)^{1-c_{i,\tau i}}[m-1]_{q_i}^!(q_i-q_i^{-1})^{m-2}q_i^{\frac{(1-m)(m-2+2c_{i,\tau i})}{2}}B_i^{(m-1)} \tk_{\tau i};
 \notag \\
\sum_{r+s=m} & (-1)^{r+c_{i,\tau i}} q_i^{-r(1-c_{i,\tau i}-m)}B_i^{(r)}B_{\tau i}B_i^{(s)} \label{eq:HOSII2}\\ \notag
&= (-1)^{m+c_{i,\tau i}+1}[m-1]_{q_i}^!(q_i-q_i^{-1})^{m-2}q_i^{\frac{(m-1)(m-2+2c_{i,\tau i})}{2}+c_{i,\tau i}}B_i^{(m-1)}\tk_i.
 \end{align}
\end{thm}

\begin{proof}
Note that
\begin{align*}
(q_i^{-2};q_i^{-2})_{-c_{i,\tau i}}=&(1-q_i^{-2})(1-q_i^{-4})\cdots (1-q_i^{2c_{i,\tau i}}) \\
	&= q_i^{\frac{-c_{i,\tau i}^2+c_{i,\tau i}}{2}}(q_i-q_i^{-1})^{-c_{i,\tau i}}[-c_{i,\tau i}]_{q_i}^!.
\end{align*}
Applying the bar involution (which sends $q^k\mapsto q^{-k}$) to the above identity, we obtain
\begin{align*}
(q_i^{2};q_i^{2})_{-c_{i,\tau i}} 
&=(-1)^{c_{i,\tau i}}q_i^{\frac{c_{i,\tau i}^2-c_{i,\tau i}}{2}}(q_i-q_i^{-1})^{-c_{i,\tau i}}[-c_{i,\tau i}]_{q_i}^!.
\end{align*}
	
So we can rewrite the BKL relation \eqref{relation5} as
\begin{align}
 \label{eq:rewr}
\sum_{r+s=1-c_{i,\tau i}} (-1)^{r+c_{i,\tau i}}&B_i^{(r)}B_{\tau i}B_i^{(s)}=[-c_{i,\tau i}]_{q_i}^!(q_i-q_i^{-1})^{-c_{i,\tau i}-1}\\
&\left(q_i^{\frac{-c_{i,\tau i}^2+3c_{i,\tau i}}{2}}B_{i}^{(-c_{i,\tau i})}\tk_i
-(-1)^{c_{i,\tau i}}q_i^{\frac{c_{i,\tau i}^2-c_{i,\tau i}}{2}}B_{i}^{(-c_{i,\tau i})}\tk_{\tau i}\right),\notag
\end{align}
which is equivalent to 
\begin{align*}
		\ty_{i,\tau i;1,1-c_{i,\tau i},e}=0.
\end{align*}
Then the first identity in \eqref{eq:y=0} follows by Theorem \ref{thm:recursion} and by induction on $m$. The second identity in \eqref{eq:y=0} follows from the first one by applying $\sigma$. 

It remains to prove \eqref{eq:HOSII1}--\eqref{eq:HOSII2}. Set $m>1-c_{i,\tau i}$ throughout the remainder of this proof. 
The identities \eqref{eq:HOSII1} and \eqref{eq:HOSII2} are equivalent by applying the bar involution $\psi_\imath$ in Lemma~\ref{lem:bar}. 

We shall show that $\ty_{i,\tau i;1,m,1}=0$ implies (and is actually equivalent to) the identity \eqref{eq:HOSII1}; similarly, $\ty_{i,\tau i;1,m,-1}=0$ is actually equivalent to the identity \eqref{eq:HOSII2}.

First note that 
\begin{align}
\label{eq:prod1}  
\prod_{j=0}^{c_{i,\tau i}+m-2} \big(-q_i^{(2j-c_{i,\tau i}-2m+2)}+q_i^{c_{i,\tau i}-2}\big)=0.
\end{align}
On the other hand, a direct computation shows that 
\begin{align}
\prod_{j=0}^{c_{i,\tau i}+m-2} & (-q_i^{(2j-c_{i,\tau i}-2m+2)}+q_i^{2-c_{i,\tau i}})
\label{eq:prod2}  \\
&= q_i^{-\frac{(c_{i,\tau i}+m-2)(c_{i,\tau i}+m-1)}{2}}(q_i-q_i^{-1})^{m+c_{i,\tau i}-1}[m]_{q_i}[m-1]_{q_i}\cdots[2-c_{i,\tau i}]_{q_i}.
\notag
\end{align}
Now the identity \eqref{eq:HOSII1} follows from $\ty_{i,\tau i;1,m,1}=0$ (see \eqref{eq:fitaui}) and \eqref{eq:prod1}--\eqref{eq:prod2}.
%
\end{proof}

The relation \eqref{relation5} and the corresponding Serre-Lusztig relation for  $c_{i,\tau i}=0$ turn out to be a variant of familiar formulas in quantum $\mathfrak{sl}_2$.

\begin{cor} [\text{cf. \cite[Corollary 3.1.9]{Lus93}}]
  \label{cor:rank1}
Let $i\in I$ be such that $c_{i,\tau i}=0$. For any $N,M\geq0$ we have in $\tUi$,
\begin{align}
       \label{eq:BB1}
B_{\tau i}^{(N)}B_i^{(M)}=\sum_{t\geq0} B_{i}^{(M-t)}\prod_{s=1}^t\frac{q_i^{2t-N-M-s+1}\tk_i-q_i^{-2t+N+M+s-1}\tk_{\tau i}}{q_i^s-q_i^{-s}}B_{\tau i}^{(N-t)};
        	\\
        	\label{eq:BB2}
B_{ i}^{(N)}B_{\tau i}^{(M)}=\sum_{t\geq0} B_{\tau i}^{(M-t)}\prod_{s=1}^t\frac{q_i^{2t-N-M-s+1}\tk_{\tau i}-q_i^{-2t+N+M+s-1}\tk_{ i}}{q_i^s-q_i^{-s}}B_{ i}^{(N-t)}.
\end{align}
\end{cor}

\begin{proof}
Thanks to $c_{i,\tau i}=0$, the relations \eqref{relation1} and \eqref{relation5} become
\begin{align*}
\tk_i B_i = q_i^{-2} B_i \tk_i, \quad
\tk_{\tau i}  B_{\tau i}  &= q_i^{-2} B_{\tau i}  \tk_{\tau i} ,
\\
\tk_{\tau i} B_i = q_i^{2} B_i \tk_{\tau i}, \quad
\tk_i B_{\tau i}  &= q_i^{2} B_{\tau i} \tk_i,
\\
B_{\tau i}B_i-B_iB_{\tau i} &=\frac{\tk_i-\tk_{\tau i}}{q_i-q_i^{-1}},
\end{align*}
That is, $\{B_{\tau i}, B_i, \tk_i, \tk_{\tau i} \}$ generate the Drinfeld double quantum group $\tU_{q_i}(\mathfrak{sl}_2)$ (where $B_i$ plays the role of $F_i$). Now the proof of \cite[Corollary 3.1.9]{Lus93} can be repeated here. 
\end{proof}

\begin{rem}
	Theorems~\ref{thm:recursion} and \ref{thm:HOSII} remain valid over $\Ui =\Ui_\bvs$, once we replace $\tk_i$ and $\tk_{\tau i}$ by $\vs_{\tau i}k_i$ and $\vs_ik_{i}^{-1}$ respectively; see Proposition \ref{prop:QSP12}.
\end{rem}

\subsection{Additional Serre-Lusztig relations}

Let $i,j \in I$ be such that $i, j$, and $ \btau i$ are all distinct. 
Recall the relation $\sum_{n=0}^{1-c_{ij}} (-1)^nB_i^{(n)}B_jB_i^{(1-c_{ij}-n)} =0$ in $\tUi$ from  \eqref{relation3}. 

The Serre-Lusztig relations associated to \eqref{relation3} take the same form as for quantum groups \cite[\S 7.1.1]{Lus93}. More explicitly, for any $m\in\N$, let
\begin{align*}
\ty_{i,j;n,m,e}:=&\sum\limits_{r+s=m}(-1)^r q_i^{er(-c_{ij}-m+1)} B_i^{(r)} B_j^{(n)} B_i^{(s)};
\\
\ty'_{i,j;n,m,e}:=&\sigma(\ty_{i,j;n,m,e}).
\end{align*}

\begin{prop}
\label{prop:LS2}
Assume that $i,j\in I$ are such that $i, \tau i,$ and $j$ are all distinct. 
Then, for $m>-c_{ij}$, we have
\begin{align*}
\ty_{i,j;n,m,e}=0,\qquad\quad  
\ty'_{i,j;n,m,e}=0.
\end{align*}
\end{prop}

\begin{proof}
By a formal computation (as in \cite[Lemma 7.1.2]{Lus93}), we have the following recursive formulas:
\begin{align*}
-q_i^{-e(2m+nc_{i,j})} & B_i \ty_{i,j;n,m,e} +\ty_{i,j;n,m,e}B_i
	=[m+1]_{i}\, \ty_{i,j;n,m+1,e},
	\\
-q_i^{-e(2m+nc_{i,j})} & \ty'_{i,j;n,m,e}B_i+B_i\ty'_{i,j;n,m,e}
	=[m+1]_{i}\, \ty'_{i,j;n,m+1,e}.
\end{align*}
By \eqref{relation3}, we have $\ty_{i,j;n,1-c_{ij},1} =0$, and thus $\ty_{i,j;n,1-c_{ij},-1} =0$ by applying the bar involution. The general formula in the proposition follows by induction on $m$ using the recursive formulas above .
\end{proof}


\begin{rem}
	The result in Proposition~ \ref{prop:LS2} also holds in $\Ui$.
\end{rem}


\subsection{Braid group symmetries for $\tUi$}

Let $W=(s_i\mid i\in I)$ be the Weyl group of $\fg$. Recall Lusztig constructed four variants of automorphisms,  $T''_{i,e}$ and $T'_{i,e}$, for $i\in I$ and $e \in \{\pm 1\}$, of the quantum group $\U$ \cite[Chapter 37]{Lus93}.

We regard $\tau$ and elements in $W$ as automorphisms on $\Z^I$. The {\em restricted Weyl group} associated to the quasi-split symmetric pair $(\fg,  \fg^{\omega \btau})$ are defined to be the following subgroup of $W$:
\begin{align}
	\label{eq:Wtau}
	W^{\btau} =\{w\in W\mid \btau w =w \btau\}.
\end{align}

Recall the subset $\I_\tau$ of $\I$ from \eqref{eq:ci}, and define
\begin{align}
	\label{eq:Itau2}
	\ov{\I}_\btau:=\{i\in\I_\btau\mid c_{i,\btau i}=-1,0, \text{ or }2 \}.
\end{align}
In our setting, $\ov{\I}_\btau$ consists of exactly those $i \in \I_\btau$ such that the $\btau$-orbit of $i$ is of finite type. We denote by $\bs_{i}$, for $i\in\ov{\I}_\btau$, the following element of order 2 in the Weyl group $W$
\begin{align} \label{eq:bsi}
	\bs_i= \left\{
	\begin{array}{ll}
		s_{i}, & \text{ if } i=\btau i
		\\
		s_is_{\btau i}, & \text{ if }c_{i,\btau i}=0,
		\\
		s_is_{\btau i}s_i, & \text{ if }c_{i,\btau i}=-1.
	\end{array}
	\right.
\end{align}
According to (a special case of) \cite[Appendix]{Lus93}, the restricted Weyl group $W^{\btau}$ can be identified with a Coxeter group with $\bs_i$ ($i\in \ov{\I}_\btau$) as its generators.

For any $i\in\bar{I}_\btau$ such that $c_{i,\tau i}=2$ (i.e., $\tau i=i$), the existence of  automorphisms $\T'_{i,e}$, $\T''_{i,e}$ on $\tUi$ together with explicit formulas for their actions on Chevalley generators are conjectured in \cite[Conjecture 6.5]{CLW21}. 

Below we make a conjecture on the existence of  automorphisms $\T'_{i,e}$, $\T''_{i,e}$ on $\tUi$ with explicit formulas, for $i\in\bar{I}_\btau$ such that $c_{i,\tau i}=0$. 

\begin{conj}	
   \label{conj:BG1}
For any $i\in\bar{I}_\btau$ such that $c_{i,\tau i}=0$,  and $e \in \{\pm 1\}$,
there are automorphisms $\T'_{i,e}, \T''_{i,e}$ on $\tUi$  such that
\begin{align*}
\T'_{i,e}(\tk_j)&=\T''_{i,e}(\tk_j)= \tk_i^{-c_{ij}} \tk_{\btau i}^{-c_{\tau i,j}} \tk_j,
\\
		\T'_{i,-1}(B_j) &=
	\begin{cases}
		-B_{\btau i}\tk_{\btau i}^{-1},  & \text{ if }j=i \\
		-\tk_{i}^{-1}B_{i},  &\text{ if }j=\btau i,
	\end{cases}\qquad\qquad
	\T'_{i,1}(B_j)=
	\begin{cases}  -B_{\btau i}\tk_{i}^{-1},  & \text{ if }j=i \\
		-\tk_{\btau i}^{-1}B_{i},  &\text{ if }j=\btau i,
	\end{cases}
\\
\T''_{i,1}(B_j) &=
\begin{cases}
 -\tk_{i}^{-1}B_{\btau i},  & \text{ if }j=i\\
-B_i\tk_{\btau i}^{-1}  ,  &\text{ if }j=\btau i,
\end{cases}
\qquad\qquad
\T''_{i,-1}(B_j)=
\begin{cases}
 -\tk_{\btau i}^{-1}B_{\btau i},  & \text{ if }j=i\\
-B_i\tk_{i}^{-1},  &\text{ if }j=\btau i,
\end{cases}
\end{align*}
and for $j\neq i,\btau i$,
\begin{align*}
	\T'_{i,-1}(B_j)
&= \sum^{-\max(c_{ij},c_{\tau i,j})}_{u=0} \; \sum^{-c_{ i,j}-u}_{r=0} \; \sum_{s=0}^{-c_{\tau i,j}-u} (-1)^{r+s} v^{r-s+(-c_{ij}-r-s-u)u } \\
&\qquad\qquad\qquad\qquad\qquad\qquad
\times  \tk_{i}^u  B_{i}^{(-c_{ij}-r-u)} B_{ \btau i}^{(s)}  B_j B_{\btau i}^{(-c_{\tau i,j}-u-s)} B_{ i}^{(r)},
\\
\T'_{i,1}(B_j)
= &\sum^{-\max(c_{ij},c_{\tau i,j})}_{u=0} \; \sum^{-c_{ i,j}-u}_{r=0} \; \sum_{s=0}^{-c_{\tau i,j}-u} (-1)^{r+s} v^{- (r-s+(-c_{ij}-r-s-u)u )} \\
&\qquad\qquad\qquad\qquad\qquad\qquad
\times  \tk_{\btau i}^u B_{i}^{(-c_{ij}-r-u)} B_{\btau i}^{(s)} B_j B_{\btau i}^{(-c_{\tau i,j}-u-s)}  B_{i}^{(r)},
\\
\T''_{i,1}(B_j)
&= \sum^{-\max(c_{ij},c_{\tau i,j})}_{u=0} \; \sum^{-c_{ i,j}-u}_{r=0} \; \sum_{s=0}^{-c_{\tau i,j}-u} (-1)^{r+s} v^{ r-s+(-c_{ij}-r-s-u)u } \\
&\qquad\qquad\qquad\qquad\qquad\qquad
 \times B_i^{(r)} B_{\tau i}^{(-c_{\tau i,j}-u-s)} B_j B_{\tau i}^{(s)} B_i^{(-c_{ij}-r-u)}\tk_{\tau i}^u,
\\
\T''_{i,-1}(B_j)
= &\sum^{-\max(c_{ij},c_{\tau i,j})}_{u=0} \; \sum^{-c_{ i,j}-u}_{r=0} \; \sum_{s=0}^{-c_{\tau i,j}-u} (-1)^{r+s} v^{- (r-s+(-c_{ij}-r-s-u)u )} \\
&\qquad\qquad\qquad\qquad\qquad\qquad
 \times B_i^{(r)} B_{\tau i}^{(-c_{\tau i,j}-u-s)} B_j B_{\tau i}^{(s)} B_i^{(-c_{ij}-r-u)}\tk_{i}^u.
\end{align*}
\end{conj}

Recall the bar involution $\psi_\imath$ and anti-involution $\sigma$ on $\tUi$ from Lemma~\ref{lem:bar}. 
The (conjectured) automorphisms $\T''_{i,e}$ and $\T'_{i,e}$ are related to each other by
\begin{align*}
\sigma\T_{i,e}' \sigma = \T_{i,-e}'', \quad
\psi_\imath \T_{i,e}' \psi_\imath = \T_{i,-e}', \quad
\psi_\imath \T_{i,e}'' \psi_\imath = \T_{i,-e}''.
\end{align*}
For a fixed $e$, the automorphisms $\T''_{i,e}$ (and respectively, $\T'_{i,e}$) are expected to satisfy the braid relations for $W^\tau$ defined in \eqref{eq:Wtau} (extending the suggestion in \cite{KP11} for $\Ui$ of finite type). 

We do not know of (conjectural) general formulas for the automorphisms in the remaining case for $c_{i,\tau i}=-1$ in \eqref{eq:Itau2}; see however \cite{KP11, WZ22} for $\Ui$ of type AIII.

\begin{rem}
\label{rem:weight}
Let us comment on the formulation of the above conjecture, focusing on $\T''_{i,1}$ only. 
Since $\bs_i =s_is_{\btau i}$ by \eqref{eq:bsi}, we know the expected weights for $\T''_{i,1}(x)$, for any $x \in \tUi$, and the $\imath$Hall algebra consideration helps to finalize the simple formulas other than $\T''_{i,1}(B_j)$. 

Recall the braid group symmetry $T_{i,1}''$ on $\U$ from \cite{Lus93}, and note that the weight of $T_{i,1}''(F_j)$ in $\U$ is equal to the weight of a corresponding Serre relation minus $\alpha_i$. By a comparison however, there is no such a clean relation between the weight for the Serre-Lusztig relation in Theorem~\ref{thm:HOSII} and the expected weight for $\T''_{i,1}(B_j)$, and hence the Serre-Lusztig relation is not helpful in our guess for formulas of $\T''_{i,1}$. (This is very different from the case treated in \cite{CLW21}.)

Using $\bs_i =s_is_{\btau i}$ again, we expect that the formula for $T_{i,1}'' T''_{\btau i,1}$ on $\U$ provides the leading term for an expected formula for $\T''_{i,1}$ on $\tUi$, and this explains the leading summands in $\T''_{i,1}(B_j)$ for $u=0$ in the above conjecture. 
Some partial case computations (e.g., $c_{ij}=0, -1$), our experience in \cite{CLW21} on ``lower terms", and some lucky  guesses allow us to pin down the conjectural general formula for $\T''_{i,1}(B_j)$. 
\end{rem}

\begin{rem}
\label{rem:weight2}
We shall develop a Hall algebra approach in \cite{LW21b} to prove Conjecture~\ref{conj:BG1} for {\em quasi-split} $\imath$quantum groups (under an additional assumption that the Cartan matrix is symmetric and $c_{j, \tau j} \in 2\Z$, for all $j$). For the Hall algebra approach to succeed, it is crucial to have our precise formulas available. 
\end{rem}

\subsection{General case}
 \label{subsec:nonQS}
 
Let us explain that the main results and conjecture are actually valid in greater generality. 
 
We recall (not necessarily quasi-split) $\imath$quantum group $\Ui$ and $\tUi$ quickly. 
Let $\tau$ be an involution of the Cartan datum $(\I, \cdot)$; we allow $\tau =\Id$.
Let $\I_{\bullet} \subset \I$ be a Cartan subdatum of {\em finite type}. Let $W_{\I_\bullet}$ be the Weyl subgroup for $(\I_{\bullet}, \cdot)$ with $w_{\bullet}$ as  its longest element. 
Denote $\I_{\circ} = \I \backslash \I_{\bullet}.$ 
The pair $(I =I_\circ \cup I_\bullet,  \tau)$ is required to satisfy some compatibility conditions (cf. \cite[Definition~2.3]{Ko14}). Associated to such an admissible pair $(I =I_\circ \cup I_\bullet, \tau)$, an $\imath$quantum group $\Ui$ is defined as a subalgebra of $\U$ \cite{Le99, Ko14}.

Following and generalizing \cite{LW19a}, we define a universal $\imath$quantum group $\widetilde{\bU}^\imath$ to be the $\K(q)$-subalgebra of the Drinfeld double $\tU$ generated by $E_\ell, F_\ell, \tK_\ell, \tK'_\ell$, for $\ell \in \Iblack$, and
\[
B_i= F_i + T_{\wb} (E_{\btau i}) \tK_i',
\qquad \tk_i= \tK_i \tK_{\btau i}', \quad \forall i \in \Iwhite.
\]
Here $T_w$, for $w\in W_{\Iblack}$, is the same as $T''_{w,+1}$ in \cite[Chapter~ 37]{Lus93}.
Then $\tUi$ is a coideal subalgebra of $\tU$. In particular, $\tUi$ contains the Drinfeld double quantum group $\tU_{I_{\bullet}}$ associated to $I_\bullet$ (generated by $E_\ell, F_\ell, \tK_\ell, \tK'_\ell$, for $\ell \in \Iblack$) as a subalgebra. 
If $I_{\bullet}=\emptyset$, then the $\imath$quantum group $\tUi$ is quasi-split as defined in   \S\ref{subsec:iQG}. 

The recursive formulas in Theorem~\ref{thm:recursion} and the Serre-Lusztig relations in Theorem~\ref{thm:HOSII} remain valid for $\imath$quantum groups $\tUi$ in this generality, under the assumption that $i\in \I_\circ$ and $\tau i \neq i =w_\bullet i$. 

We expect Conjecture \ref{conj:BG1} (where $i, j\in I_\circ$) to be valid for a general $\imath$quantum group $\tUi$ under the assumptions that $c_{i,\tau i}=0$ and $i =w_\bullet i$; in addition, $\T''_{i,e}$ and $\T'_{i,e}$ fix (the generators of) the subalgebra $\tU_{I_{\bullet}}$. 
It is worth noting that the braid group actions on $\imath$quantum groups of all finite types are now constructed in \cite{WZ22}, confirming the main conjecture in \cite{KP11}.


\end{document}